\documentclass[10pt]{article}
\usepackage{amsmath}
\usepackage{amsfonts}
\usepackage{amssymb}
\usepackage{paralist}
\usepackage{graphics} %% add this and next lines if pictures should be in esp format
\usepackage{epsfig} %For pictures: screened artwork should be set up with an 85 or 100 line screen
\usepackage[colorlinks=true]{hyperref}

\hypersetup{urlcolor=blue, citecolor=red}

\newtheorem{theorem}{Theorem}[section]
\newtheorem{corollary}{Corollary}[section]

\newtheorem{lemma}{Lemma}[section]
\newtheorem{proposition}{Proposition}[section]

\newtheorem{definition}{Definition}[section]
\newtheorem{remark}{Remark}[section]

\newcommand{\ep}{\varepsilon}

\def\sqr#1#2{{\vcenter{\vbox{\hrule height.#2pt
              \hbox{\vrule width.#2pt height#1pt \kern#1pt \vrule
width.#2pt}
              \hrule height.#2pt}}}}
\def\3n{\negthinspace \negthinspace \negthinspace }
\def\2n{\negthinspace \negthinspace }
\def\1n{\negthinspace }

\def\liminf{\mathop{\underline{\rm lim}}}

\def\({\Big (}
\def\){\Big )}
\def\[{\Big[}
\def\]{\Big]}

\makeatletter
   
   \@addtoreset{equation}{section}
\makeatother

\def\ba{\begin{array}}
\def\be{\begin{equation}}
\def\bc{\begin{corollary}}
\def\bd{\begin{definition}}
\def\bea{\begin{eqnarray}}
\def\bel{\begin{equation}\label}
\def\bl{\begin{lemma}}
\def\bp{\begin{proposition}}
\def\br{\begin{remark}}
\def\bt{\begin{theorem}}
\def\ee{\end{equation}}
\def\eea{\end{eqnarray}}
\def\et{\end{theorem}}
\def\ec{\end{corollary}}
\def\el{\end{lemma}}
\def\ep{\end{proposition}}
\def\er{\end{remark}}
\def\ea{\end{array}}
\def\ed{\end{definition}}

\title%[LQ Time-Inconsistent Optimal Control]
      {Consensus Seeking in Multi-Agent Systems with  Multiplicative Measurement Noises}

\author{Yuan-Hua Ni\thanks{Department of Mathematics, Tianjin Polytechnic University, Tianjin, P.R. China ({\tt yhni@amss.ac.cn}).}~~~~~~ Xun Li\thanks{Department of Applied Mathematics, The Hong Kong Polytechnic University, Hunghom, Kowloon, Hong Kong (malixun@polyu.edu.hk).}}

\begin{document}
\maketitle

%The abstract of your paper
\begin{abstract}

In this paper, the consensus problems of the continuous-time
integrator systems under noisy measurements are considered. The measurement noises, which appear when agents measure their neighbors' states, are modeled to be multiplicative. By multiplication of
the noises, here, the noise intensities are proportional to the absolute value of the relative states of agent and its neighbor. By using known
distributed protocols for integrator agent systems, the closed-loop {system
is} described in the vector form by a singular stochastic differential equation. For the fixed and switching network topologies cases,
constant consensus gains are properly selected, such that mean square
consensus and strong consensus can be achieved. Especially, exponential mean square convergence of agents' states to the common value is derived for the fixed topology case. In addition, asymptotic
unbiased mean square average consensus and asymptotic unbiased
strong average consensus are also studied.
Simulations shed light on the effectiveness of the proposed
theoretical results.

\end{abstract}

%The title of your section 1
\section{Introduction}{\label 1}

Recently, distributed coordination of multi-agent systems has
attracted more and more attention of multi-disciplinary researchers,
due to its wide applications in cooperative control, formation
control and distributed optimization, and flocking problem. Most of
all, the consensus problem is one of the fundamental topics in
distributed coordination. By consensus, the group of
dynamic agents will asymptotically reach an agreement on certain quantity of interest. For the distributed control problem,
this means that by designing a distributed protocol such that the states of  each
agent asymptotically reach the agreement. This is the core topic of this
paper.

The research efforts of consensus problem in the system and control community can be traced back to the work \cite{Tsitsiklis}, which deals with the asynchronous
consensus problems with application to distributed decision-making
systems. Much recent works in this area are motivated by the Vicsek's model \cite{Vicsek}. In \cite{Vicsek}, Vicsec \emph{et al}. propose a nearest neighbor rule that updates the heading of the autonomous agents moving in the plane
with the same speed but with different headings, and provide simulation results which demonstrate that
the nearest neighbor rule can cause all agents
to eventually move in the same direction. For the Vicsek model, its analytic behavior is subsequently studied in \cite{Jadbabaie}, which provides a theoretical explanation for the observed
behavior. Another exact formulation in population of autonomous agents is achieved by \cite{Smale}. The cardinal feature of the model in \cite{Smale} is that the interaction between any pair of the agents is a well defined nonlinear function of the difference between their co-ordinates in $\mathbb{R}^3$.  Since then, lots of literatures about this area appear,  and readers may refer to, for example, \cite{Aysal,Moreau,Ni-wei,Murray,Ren wei-0,Ren Wei,Wang-long}, and to \cite{Olfati,Ren wei-2,Vicsek-2}  for recent survey.

In reality, communication processes are always corrupted by various
uncertain factors, such as random link failures, transmission noises
and quantization errors.
These are outcomes of the use of sensors, quantization and
wireless fading channels in the network. Recently, the consensus
problems that concern corrupted communication processes between agents have attracted many
researchers. In \cite{Huang-09-SIAM}, the authors model the measurement noises to be additive, which means that the noises additively input the communication processes.
To attenuate this type of measurement noises, a decreasing consensus gain is designed to reduce the detrimental effect of the noises. For fixed  network topology, mean square consensus and strong consensus results
are presented. Furthermore, in \cite{Huang-10-Auto}, the authors deal with randomly varying
topology, while \cite{Huang-10-Automatica} considers the Markovian
and arbitrary switching topologies. In \cite{Li-Tao-2009-Automatica}, the authors consider continuous time first-order integrator model with Gaussian additive standard white noise. To be exact, in \cite{Li-Tao-2009-Automatica}, firstly, for noise-free
cases, necessary and sufficient conditions are given on the network topology and consensus gains to
achieve average-consensus; secondly, for the cases with measurement noises, necessary and sufficient
conditions are given on the consensus gains to achieve asymptotic unbiased mean square average-consensus. In addition, \cite{Li-Tao-2010--IEEE} considers the time-varying topology case, while \cite{Ma-cui-qin}
deals with the leader-follower consensus control problem. {It is worth noting that most existing literatures are for additive measurement noises.}

In this paper, we model the measurement noises {to be multiplicative}, which may be viewed as the complement to those considered in \cite{Huang-09-SIAM,Li-Tao-2009-Automatica,Ma-cui-qin}. {Our modelling comes from a simple intuition. To be exact}, as the states of agents (say, for example, mobile vehicles) may be viewed as the positions, intuitively, the closer the vehicles are to each other, the smaller the intensities of measurement noises should be related to them. {In fact, the idea that takes into account the relative states has been already proposed by Cucker and Smale in \cite{Smale}, which is called an exact formulation by Vicsek \cite{Vicsek-2}.} In \cite{Smale}, elements of the adjacency matrix are modeled as non-increasing functions of the distance between agents to measure the ability of agents communication with each other. This means that the closer the agents are to each other, the bigger the influences is among them.
{Furthermore, many interesting extensions of Cucker-Smale model are developed; see, for example, \cite{Carrillo,Cucker,Ha,Vicsek-2}. On the other hand, multiplicative noises often appear in communication process modelling, especially for the cases of fading communication channel; see, for example, \cite{Dogandzic,Jerry,Rong,Willems}. Therefore, in this paper, to describe the influence of relative states on the communication ability between agents under the noisy environment, we model the intensities of measurement noises to be proportional to the distance  between agents. Clearly, our model may be viewed as a variant of Cucker-Smale model in the case of noisy environment, and the noises are multiplicative in the way of depending on the agents' relative states. Furthermore, the measurement noises will disappear when states of two agents coincide; in other words, agents' states can be measured precisely by agents themselves, which is the case of most existing literatures. }

In this paper, we consider consensus problems for continuous-time
first-order  integrator system. The major difference between this paper and existing
literatures is that  the measurement model  is different. This results in that  the closed-loop
system is a singular stochastic differential equation (SDE) under the known protocol. Here, by
singularity of SDE, we mean that the covariance matrix of the system noise
is allowed to be degenerate. This is different from the non-degenerate SDE arising in existing literatures due to the additivity of the measurement
noises. Intuitively, if
we can design a mean square (strong, respectively) protocol, and when time $t$ is
large enough, the scales of measurement noises will be very small. Therefore, the
decreasing consensus gain is not necessary. In fact, methodology developed here is related to the stability of equilibrium point 0 of linear stochastic systems with multiplicative noises \cite{zhou,Damm,zhangweihai,zhangweihai-3}.
This paper considers two cases
of network topology. For the
case of fixed  network topology, a constant consensus gain is
 selected to ensure the mean square consensus and strong
consensus, respectively. Especially, exponential mean square convergence of agents' states to the common value is derived. For the case of switching network topology, a consensus gain is also designed. If all the possible digraphs are balanced and the union of them contains a spanning tree, we prove by contradiction that the expected total consensus error will approach to zero asymptotically. Therefore, the mean square consensus and the strong consensus are achieved easily.
%On the other hand, we prove similar results under an alternative condition such that the balancedness of the possible %digraphs is dropped away.

The remainder of this paper is organized as follows. Section
\ref{problem formulation} contains formulation of the problem. Consensus problem under fixed topology is discussed in Section \ref{section-fixed}. Section
\ref{section-switching} generalized the results obtained in Section
\ref{section-fixed} to the case of switching topology. Illustrative examples are presented in Section
\ref{simulation}, and Section \ref{conclusion} gives some concluding
remarks.

\section{Problem formulation}\label{problem formulation}

Consider a multi-agent system consisting of $N$ agents,
labeled by agents $1,2,\cdots,N$. Each agent's dynamics is the
continuous-time first-order integrator
\begin{eqnarray}\label{agent-dynamics}
\dot{x}_i(t)=u_i(t),~~i=1,2,...,N,
\end{eqnarray}
where $x_i(t), u_i(t)\in \mathbb{R}$,  are the state  and control input of agent $i$.
$u_i(t)$ is designed only based on the local information of neighbors of agent $i$. Here, agent $j$ is called a neighbor of agent $i$, if
agent $i$ can receive measurement of state of agent $j$.
If agent $j$ is a neighbor of agent $i$, we may denote this relationship by an ordered pair $(j,i)$, called an edge from agent $j$ to agent $i$. Furthermore, if the agents are viewed as abstract nodes, the information exchanging among agents at time $t$ may be modeled as a
directed graph (digraph) $G(t)=(\mathcal{N},\mathcal{E}(t))$. Here, $\mathcal{N}=\{1,2,...,N\}$, while $\mathcal{E}(t)$ is the totality of all the edges at time $t$. $G(t)$ is called the network topology (or topology) of the multi-agent system at time $t$, while $x(t)\equiv [x_1(t),\cdots,x_N(t)]^T$ is called the information state of $G(t)$ at time $t$. The totality of all the neighbors of agent $i$ at time $t$ is denoted by $\mathcal{N}_{ti}$. For more about of graph theory, readers may refer to standard textbooks, or to \cite{Huang-09-SIAM,Li-Tao-2009-Automatica} that relate to the topic of this paper.

Let agent $j$ be a neighbor of agent $i$. In this paper, we model the measurement of state of agent $j$ received by agent
$i$ by
\begin{eqnarray}\label{measurement-0}
y_{ji}(t)=x_j(t)+\sigma_{ji}|x_j(t)-x_i(t)|\xi_{ji}(t),
\end{eqnarray}
where $\xi_{ji}\equiv \{\xi_{ji}(t), t\geq 0\}$, $i=1,2,...,N, j\in\mathcal{N}_{ti}$, are standard white
noises, $\sigma_{ji}\geq 0$. By (\ref{measurement-0}), when $j=i$,
$y_{ii}(t)=x_i(t)$, which indicates that agent $i$ can measure its
state $x_i(t)$ exactly.
A group of controls
$\{u_i,i=1,2,...,N\}\equiv u$ is called an admissible distributed
protocol, if $u_i(t)$ is measurable with respect to the $\sigma$-algebra
generated by $(x_i(s),y_{ji}(s),j\in \mathcal{N}_{si},0\leq s\leq t)$,
i.e., $u_i(t)\in \sigma(x_i(s),y_{ji}(s),j\in \mathcal{N}_{si},0\leq
s\leq t)$. The so-called consensus seeking problem is to design an admissible distributed protocol, such that the states of all the agents asymptotically approach a common value in some sense. If the convergence is in the mean square sense (almost sure sense, respectively), we call the multi-agent system achieves mean square (strong) consensus, and the corresponding protocol is called a mean square (strong) distributed protocol.

For (\ref{agent-dynamics}), a known protocol is
\begin{eqnarray}\label{protoco0}
u_i(t)=a\sum_{j\in \mathcal{N}_{ti}}\left(y_{ji}(t)-x_i(t)\right),
\end{eqnarray}
$i=1,2,...,N$, which is very popular in the noise-free case (i.e., $y_{ji}(t)=x_j(t)$). In (\ref{protoco0}), if $\mathcal{N}_{ti}=\varnothing$, $u_i(t)$ is viewed as $0$. In this paper, protocol (\ref{protoco0}) with $y_{ji}(t)$ given by (\ref{measurement-0}) is also taken. Here, $a$ is called the consensus gain, which will be determined below. To ease the following analysis, the adjacency matrix $\mathcal{A}(t)$ of $G(t)$ is introduced: for any $i,j\in\mathcal{N}$, the element {$a_{ji}(t)$} of $\mathcal{A}(t)$ is equal to 1 if and only if $(j,i)\in\mathcal{N}_{ti}$, otherwise, {$a_{ji}(t)=0$.} By (\ref{agent-dynamics})(\ref{measurement-0})(\ref{protoco0}), we have
\begin{eqnarray}\label{closed-loop-0-0}
\dot{x}_i(t)=a\sum_{j\in\mathcal{N}}a_{ji}(t)(x_j(t)-x_i(t))+a\sum_{j\in\mathcal{N}}a_{ji}(t)\sigma_{ji}|x_j(t)-x_i(t)|\xi_{ji}(t),~~i=1,2,...,N,
\end{eqnarray}
where the fact that $a_{ji}(t)=0$ if $(j,i)\not\in\mathcal{N}_{ti}$ is used. Now, define a matrix $\Sigma(t)\in \mathbb{R}^{N\times N^2}$, whose $i$-th line is given by $(0,..,0,a_{1i}(t)\sigma_{1i},a_{2i}(t)\sigma_{2i},...,a_{Ni}(t)\sigma_{Ni},0,...0)$ with $a_{1i}(t)\sigma_{1i}$ being the $((i-1)N+1)$-th element. Let $y_i(t)=\mbox{diag}(|x_1(t)-x_i(t)|,|x_2(t)-x_i(t)|,...,|x_N(t)-x_i(t)|)$, $y(t)=\mbox{diag}(y_1(t),y_2(t),...,y_N(t))$.
Denote $\xi_i(t)=[\xi_{1i}(t),\xi_{2i}(t),...,\xi_{Ni}(t)]^T$, and $\xi(t)=[(\xi_{1}(t))^T, (\xi_{2}(t))^T,$ $...,(\xi_{N}(t))^T]^T$. The vector form of (\ref{closed-loop-0-0}) is
\begin{eqnarray}\label{closed-loop-0-1}
\dot{x}(t)=-aL(t)x(t)+a\Sigma(t) y(t)\xi(t),
\end{eqnarray}
where $L(t)$ is the Laplacian matrix of $G(t)$. This is a random differential equation. By classic stochastic analysis theory,
it makes sense to consider the integral form of (\ref{closed-loop-0-1}) in the It\^o sense
\begin{eqnarray*}
x(t)=x(0)-\int_0^taL(s)x(s)ds+\int_0^ta\Sigma(s) y(s)dW(s),
\end{eqnarray*}
i.e.,
\begin{eqnarray}\label{closed-loop-0-2}
dx(t)=-aL(t)x(t)dt+a\Sigma(t) y(t)dW(t),
\end{eqnarray}
where $\{W(t),t\geq 0\}$ is the $\mathbb{R}^{N^2}$-valued standard Brownian motion defined on a probability  space $(\Omega,\mathcal{F},P)$.  Notice that the covariance of the noise term is $a^2\Sigma(t) y(t)(\Sigma(t) y(t))^T$, which may be degenerate. Therefore, (\ref{closed-loop-0-2}) may be a singular stochastic differential equation.

\section{Networks with Fixed Topology}\label{section-fixed}

In this section, we deal with the case that the network topology $G(t)$ is fixed, i.e., $G(t)\equiv G$, for some digraph $G=(\mathcal{N}, \mathcal{E})$. Therefore, in this case, the closed-loop system is
\begin{eqnarray}\label{closed-loop-1}
dx(t)=-aLx(t)dt+a\Sigma y(t)dW(t).
\end{eqnarray}
To this end, we need the following assumption.

(\textbf{A1}). The digraph $G$ contains a spanning tree.

Note that $-L$ may be interpreted as the generator of a continuous time Markov chain with state space being $\mathcal{N}$. Therefore, some standard results on Markov chain can be used to simplified our analysis. To be exact,
by results in \cite{Huang-10-Auto}, there exists a nonsingular matrix
$\Phi=(1_N, \Phi_2)$
such that
\begin{eqnarray}\label{L}\Phi^{-1}(-L)\Phi= \left(\begin{array}{cc}0&0 \\0& \widetilde{L}\end{array}\right),\end{eqnarray}
and $\widetilde{L}\in \mathbb{R}^{(N-1)\times (N-1)}$ is Hurwitz,
where $1_N$ is a column vector with all $N$ entries being $1$,
$\Phi_2$ is a  $N\times (N-1)$ matrix. Clearly, there exists a positive definite matrix $Q$ such that
\begin{eqnarray}Q\widetilde{L}+\widetilde{L}^TQ=-I_{N-1}.
\end{eqnarray}
In addition, $\Phi^{-1}$ has
the following form
\begin{eqnarray}\label{inverse-Phi}
\Phi^{-1}=\left(\begin{array}{c}\pi\\\Psi_2
\end{array}\right),
\end{eqnarray}
where $\Psi_2$ is a $(N-1)\times N$ matrix and $\pi$ is the unique
invariant probability measure of the Markov chain with respect to the generator $-L$.
Let $z(t)=\Phi^{-1}x(t)\equiv [z^1(t), \widetilde{z}(t)^T]^T$, where $z^1(t)\in \mathbb{R}, \widetilde{z}(t)\in \mathbb{R}^{N-1}$. Thus
\begin{eqnarray}\label{x(t)}
x(t)=\Phi z(t)=z^1(t)1_N+\Phi_2\widetilde{z}(t).\end{eqnarray}
Equivalently,
$x_i(t)=z^1(t)+\phi_i\widetilde{z}(t),$
where $\phi_i$ is the $i$th row of $\Phi_2$. Clearly,
$|x_i(t)-x_j(t)|=|(\phi_i-\phi_j)\widetilde{z}(t)|. $
Therefore, $y(t)$ is a simple function of $\widetilde{z}(t)$. It is worth noting that $z^1(t)$ does not appear in this equality.
Upon the above notations, (\ref{closed-loop-1}) is equivalent to
\begin{eqnarray}\label{equation-1}
\left\{\begin{array}{l}dz^1(t)=a\beta(\widetilde{z}(t))dW(t),\\
\label{equation-2}d\widetilde{z}(t)=a\widetilde{L}\widetilde{z}(t)dt+a\widetilde{\Sigma}(\widetilde{z}(t))dW(t),
\end{array}\right.\end{eqnarray}
where
$\Phi^{-1}\Sigma y(t) =\left(\begin{array}{c}\pi \Sigma y(t)\\ \Psi_2 \Sigma y(t)\end{array}\right) \equiv \left(\begin{array}{c}
\beta(\widetilde{z}(t))\\ \widetilde{\Sigma}(\widetilde{z}(t))\end{array}\right)
$. By known results about the estimation of the solution to SDE, see, for example, \cite{Lipster-0}, for any ${t}>0$ and $m\geq 1$,  there exists a constant $c_{t,m}$, such that
\begin{eqnarray}\label{martingle-1}
E|\widetilde{z}(s)|^{2m}\leq c_{t,m}(1+E|\widetilde{z}(0)|^{2m})e^{c_{t,m}t}, ~~0\leq s\leq {t}.
\end{eqnarray}

{Before stating the main result of this section, we recall the following notion (\cite{Kolmogrov}\cite{Xia-daohang}).
\begin{definition}
Let $[a, b]$ {be} a subinterval of $ \mathbb{R}$. A function $f: [a, b]\mapsto R$ is called absolutely continuous, if for any $\varepsilon>0$, there exists a $\delta>0$ such that for any finite sequence of pairwise disjoint subintervals $(a_k, b_k)$ of $[a, b]$ satisfies
$\sum_{k}(b_k-a_k)<\delta$
then
\begin{eqnarray*}
\sum_{k}|f(b_k)-f(a_k)|<\varepsilon.
\end{eqnarray*}
\end{definition}
}
{For any absolutely continuous functions $f$, it has bounded variation, and thus is differentiable almost everywhere (a.e.) with respect to Lebesgue measure.  In addition to this, we also have the following two results (\cite{Kolmogrov}\cite{Xia-daohang}).
}
{
\begin{lemma}\label{Lemma-differentiable-1}
Let $f(t)$ be a absolutely continuous function on $[a, b]$, and $\frac{df(t)}{dt}=0$, a.e., then $f(t)\equiv c$ for some $c$.
\end{lemma}
\begin{lemma}\label{Lemma-differentiable-2}
The Newton-Leibniz formula
\begin{eqnarray*}
f(t_2)-f(t_1)=\int_{t_1}^{t_2}g(s)ds
\end{eqnarray*}
works if and only if f(t) is absolutely continuous, and $g(t)=\frac{df(t)}{dt}$, a.e..
\end{lemma}
}

{Note that Lemma \ref{Lemma-differentiable-2} is a generalization of the classical Newton-Leibniz formula when $g$ is continuous.}

\begin{theorem} \label{theorem-1}
Let $\bar{a}=\frac{1}{\sum_{i=1}^N[\Psi_2^T Q\Psi_2]_{ii}\sum_{j\in\mathcal{N}_i} \left(\sigma_{ji}|\phi_j-\phi_i|\right)^2}$ with $[\Psi_2^T Q\Psi_2]_{ii}$ being the $(i,i)$-th element of matrix $\Psi_2^T Q\Psi_2$. Under assumption (A1), if the consensus gain $a$ is selected such that $0<a<\bar{a}$, then distributed protocol (\ref{protoco0}) is a mean square protocol, and the convergence of states of agents to the common value is exponential with rate $\frac{a\bar{a}-a^2}{\bar{a}\lambda_{max}(Q)}$, where $\lambda_{max}(Q)$ is the maximum eigenvalue of $Q$.
\end{theorem}

\emph{{Proof}}.  Let
$V(t)=\widetilde{z}^T(t)Q\widetilde{z}(t)$. By It\^o's formula, we have
\begin{eqnarray}\label{thm-0}
dV(t)=\left[-a|\widetilde{z}(t)|^2+{a^2}Tr\left[Q\widetilde{\Sigma}
(\widetilde{z}(t))\left(\widetilde{\Sigma}(\widetilde{z}(t))\right)^T
\right]\right]dt+2a\widetilde{z}(t)^TQ\widetilde{\Sigma}(\widetilde{z}(t))dW(t),
\end{eqnarray}
{By (\ref{martingle-1}), we know that $\int_0^t2a\widetilde{z}(s)^TQ\widetilde{\Sigma}(\widetilde{z}(s))dW(s)$ is a martingale. Therefore, we have
\begin{eqnarray}\label{thm-0-1}
EV(t)=V(0)+\int_0^tE\left[-a|\widetilde{z}(s)|^2+{a^2}Tr\left[Q\widetilde{\Sigma}
(\widetilde{z}(s))\left(\widetilde{\Sigma}(\widetilde{z}(s))\right)^T
\right]\right]ds.
\end{eqnarray}
}

{
From Lemma \ref{Lemma-differentiable-2}, we know that $EV(t)$ is absolutely continuous, and for a.e. $t\geq 0$
\begin{eqnarray}\label{thm-6}
\frac{dEV(t)}{dt}=E\left[-a|\widetilde{z}(t)|^2+{a^2}Tr\left[Q\widetilde{\Sigma}
(\widetilde{z}(t))\left(\widetilde{\Sigma}(\widetilde{z}(t))\right)^T
\right]\right].
\end{eqnarray}
}
Notice that
\begin{eqnarray*}
(\Sigma y(t))(\Sigma y(t))^T=\mbox{diag}\left(\sum_{j=1}^N\left(a_{j1}(t)\sigma_{j1}|x_j(t)-x_1(t)|\right)^2,...,\sum_{j=1}^N\left(a_{jN}(t)\sigma_{jN}|x_j(t)-x_N(t)|\right)^2\right),
\end{eqnarray*}
which is a diagonal matrix. Simple calculation shows that

\begin{eqnarray}\label{tr}
&&Tr\left[Q\widetilde{\Sigma}(\widetilde{z}(t))\left[ \widetilde{\Sigma}(\widetilde{z}(t))\right] ^T\right]
=Tr\left[\Psi_2^T Q\Psi_2(\Sigma y(t))(\Sigma y(t))^T \right]\nonumber \\
&&=\sum_{i=1}^N[\Psi_2^T Q\Psi_2]_{ii}\sum_{j\in\mathcal{N}_i} \left(\sigma_{ji}|x_j(t)-x_i(t)|\right)^2\nonumber\\
&&=\sum_{i=1}^N[\Psi_2^T Q\Psi_2]_{ii}\sum_{j\in\mathcal{N}_i} \left(\sigma_{ji}|(\phi_j-\phi_i)\widetilde{z}(t)|\right)^2\nonumber\\
&&\leq \left(\sum_{i=1}^N[\Psi_2^T Q\Psi_2]_{ii}\sum_{j\in\mathcal{N}_i} \left(\sigma_{ji}|\phi_j-\phi_i|\right)^2\right)|\widetilde{z}(t)|^2\equiv c_1|\widetilde{z}(t)|^2.\end{eqnarray}
Therefore, if the consensus gain $a$ is selected such that
\begin{eqnarray*}\gamma_1 \equiv {a}-c_1 {a}^2>0,
\end{eqnarray*}
equivalently,
\begin{eqnarray}\label{thm-4}
0<{a}<\frac{1}{c_1}=\bar{a},
\end{eqnarray}
{then we have for a.e. $t\geq 0$
\begin{eqnarray}\label{thm-8}
\frac{dEV(t)}{dt}\leq -\gamma_1EV(t).
\end{eqnarray}
}
Clearly, we have
\begin{eqnarray}\label{thm-7}
\frac{d\left(e^{\gamma_1 t}EV(t)\right)}{dt}=e^{\gamma_1 t}\left(\frac{dEV(t)}{dt}+\gamma_1 EV(t)\right)\leq 0,~~ a.e..
\end{eqnarray}
{From Lemma \ref{Lemma-differentiable-1}, we know that the solution $EV(t)$ satisfies
$e^{\gamma_1 t}EV(t)\leq V(0)$ for all $t\geq 0$, i.e.,}
\begin{eqnarray*}\label{V-inequality}
EV(t)\leq V(0)e^{-\gamma_1 t},~~t\geq 0.
\end{eqnarray*}
{Therefore,}
\begin{eqnarray}\label{thm-2}
E|\widetilde{z}(t)|^2\leq c_2 e^{-\gamma_1 t},\end{eqnarray}
for some $c_2>0$. Clearly, $\lim_{t\rightarrow\infty}E|\widetilde{z}(t)|^2=0$.

On the other hand, by (\ref{equation-1}) ,
\begin{eqnarray}\label{thm-3}
z^1(t)=z^1(0)+\int_0^ta\beta(\widetilde{z}(s))dW(s).\end{eqnarray}
By (\ref{thm-2})(\ref{thm-3}), it follows that
\begin{eqnarray}\label{thm-5}
\sup_{t\geq 0}E|z^1(t)|^2&\leq& 2\sup_{t\geq 0}\left(E|z^1(0)|^2+E|\int_0^ta\beta(\widetilde{z}(s))dW(s)|^2\right)\nonumber\\
&\leq& 2E|z^1(0)|^2+2\bar{c}_1\int_0^\infty E|\widetilde{z}(s)|^2dt<\infty,
\end{eqnarray}
for some positive constant $\bar{c}_1$. By Lyapunov inequality, $\sup_{t\geq 0}E|z^1(t)|$ must be bounded. Due to martingale convergence theorem \cite{Lipster-0}, it follows that as $t\rightarrow \infty$, $z^1(t)$ will converge almost surely to
\begin{eqnarray*}\label{z*}
z^*\equiv z^1(0)+\int_0^\infty a\beta(\widetilde{z}(s))dW(s)=\pi x(0)+\int_0^\infty a\beta(\widetilde{z}(s))dW(s).\end{eqnarray*}
In addition, by (\ref{thm-2}), we have there exists a positive number $\bar{\bar{c}}_1$ such that
\begin{eqnarray}\label{z^1(t)}
E|z^1(t)-z^*|^2\leq \bar{c}_1\int_t^\infty E|\widetilde{z}(s)|^2ds\leq \bar{\bar{c}}_1e^{-\gamma_1t},
\end{eqnarray}
which implies  $z^1(t)$ converges exponentially to $z^*$ in mean square sense. Therefore, by (\ref{x(t)})(\ref{thm-2})(\ref{z^1(t)}),  there exists a positive number ${c}_3$ such that
\begin{eqnarray*}
E|x(t)-z^*1_N|^2&\leq&
 2(E|z^1(t)1_N-z^*1_N|^2+E|\Phi_2\widetilde{z}(t)|^2)\leq {c}_3e^{-\gamma_1t}.
\end{eqnarray*}
Therefore, we conclude that the distributed protocol (\ref{protoco0}) is a mean square consensus protocol, and the convergence of states of agents to the common value is exponential with rate $\gamma_1=\frac{a\bar{a}-a^2}{\bar{a}\lambda_{max}(Q)}$.  This completes the proof.  \hfill $\Box$

To derive the strong consensus result, we need  the following lemma. Firstly, Denote
\begin{eqnarray*}
\{Z\rightarrow \}=\{\omega\in \Omega:~ \lim_{t\rightarrow\infty}Z(t,
\omega)~ \mbox{exists~and~is~finite}\},
 \end{eqnarray*}
and $Z(\infty)=\lim_{t\rightarrow\infty}Z(t, \omega)$, $\omega \in
\{Z\rightarrow \}$. The following lemma can be found in \cite{Lipster}.
\begin{lemma}\label{Lemma-------lipster}
Let $A_1$ and $A_2$ be nondecreasing processes, and let $Z$ be
a nonnegative semimartingale with $E(Z)<\infty$  and
\begin{eqnarray*}Z(t)=Z(0)+A_1 (t)-A_2(t)+M(t),~~ t\geq 0, \end{eqnarray*}
where $M$ is a local martingale. Then
\begin{eqnarray*}
\{\omega:~A_1(\infty)<\infty\}\subseteq \{Z\rightarrow \}\cap
 \{\omega:~A_2(\infty)<\infty\}.
 \end{eqnarray*}

\end{lemma}

{By Theorem \ref{theorem-1}, we know that $EV(t)$ will approach to 0 asymptotically, while Lemma \ref{Lemma-------lipster}  says that $V(t)$ converges, a.s., to a random variable. By the subsequence method, we may prove that $V(t)$ converges, a.s., to 0, and thus the strong consensus achieves. Therefore, we have the following theorem. }

\begin{theorem}\label{theorem-2}
Under the conditions of Theorem \ref{theorem-1}, the distributed
protocol (\ref{protoco0}) is a strong consensus protocol.
\end{theorem}

\emph{{Proof}}. Firstly, by Chebyshev's inequality,
it follows that $P(|V(t)-0|>\varepsilon)=P(V(t)>\varepsilon)\leq \frac{EV(t)}{\varepsilon}$. As $\lim_{t\rightarrow \infty}EV(t)=0$, we have that $V(t)$ converges to 0 in probability. Therefore, there exists a sequence of time $\{t_n,n=1,2,3,...\}$ such that $V_{t_n}$ converges to 0 almost surely.
On the other hand, by (\ref{thm-0}), (\ref{thm-4}) and Lemma \ref{Lemma-------lipster}, we assert that there exists a random variable $V^*\geq 0$ such that $\lim_{t\rightarrow\infty}V(t)=V^*, $ a.s.. By these facts, we can conclude that $V^*$ must equal to $0$, a.s., which implies that $\lim_{t\rightarrow \infty}\widetilde{z}(t)=0$, a.s..
By the analysis of Theorem \ref{theorem-1}, $z^1(t)$ converges to
$z^*$ a.s.. Therefore, by (\ref{x(t)}), we can conclude that protocol (\ref{protoco0}) is a strong consensus protocol. This completes the proof. \hfill$\Box$

Now, let us consider the average consensus problem. The following definition can be found in \cite{Li-Tao-2009-Automatica}.

\begin{definition}\label{Asymptotic unbiased mean square average consensus protocol} A distributed protocol $u$ is called an asymptotic unbiased mean square (strong, respectively) average consensus protocol if  this protocol is a mean square (strong) consensus protocol, and in addition, the corresponding group decision value $x^*$ satisfies the following properties: $Ex^*=\frac{1}{N}\sum_{i=1}^Nx_i(0)$, $Var(x^*)<\infty$.
\end{definition}

{By Theorem \ref{theorem-1} and Theorem \ref{theorem-2}, we need only to show $z^*$ satisfies $Ez^*=\frac{1}{N}\sum_{i=1}^Nx_i(0)$ and  $Var(z^*)<\infty$. The following result is clear. }

\begin{corollary}\label{corollary-1}
 Under the conditions of Theorem \ref{theorem-1} and that the digraph $G$ is balanced, then  the distributed
protocol (\ref{protoco0}) is an asymptotic unbiased mean square (strong, respectively) average consensus protocol.
\end{corollary}

\emph{{Proof}}.  By (\ref{z*}), we have that the expectation of $z^*$ is
\begin{eqnarray}\label{expection}
Ez^*=\pi x(0)=\sum_{i=1}^N\pi_ix_i(0),
\end{eqnarray}
and the variance of $z^*$ is
\begin{eqnarray*}
Var(z^*)= E\left(\int_0^\infty a\beta(\widetilde{z}(s))dW(s)\right)^2 <\infty.
\end{eqnarray*}
Therefore, by definition of asymptotic unbiased mean square (strong, respectively) average consensus,  we need only to validate that $\pi=\frac{1}{N}1_N^T$. Notice that $\pi L=0$ and $\sum_{i=1}^N\pi_i=1$. Clearly,  by Theorem 6 of \cite{Murray}, $1^T_NL=0$ is equivalent to digraph $G$  is balanced. Therefore,   if  digraph $G$ is balanced and by the uniqueness of the invariant probability measure of the Markov chain associated with generator $-L$, we must have that $\pi=\frac{1}{N}1^T_N$. This completes the proof.  \hfill $\square$

The balancedness of  digraph $G$ is quite standard in deterministic
average consensus problem, see for example \cite{Murray}. On the
other hand, we notice that the above mentioned consensus properties are all global notions, as the initial values of states of agents may vary in the whole space of  $\mathbb{R}^N$. In some sense, this suggests that the
balancedness property is necessary to obtain the ``global" average
consensus results. In fact, it is easy to show that ``local" average
consensus results may also be derived even if the balancedness
property is not satisfied. It is interesting that this phenomenan
has been hardly discussed even in literatures about deterministic
average consensus problem. We give a simple description about this here. Clearly, to achieve the average consensus, by
(\ref{expection}), the following is necessary
\begin{eqnarray*}
\sum_{i=1}^N\pi_ix_i(0)=\frac{1}{N}\sum_{i=1}^Nx_i(0).
\end{eqnarray*}
Equivalently,
\begin{eqnarray}\label{pi}
(\pi_1-\frac{1}{N})x_1(0)+\cdots+(\pi_N-\frac{1}{N})x_N(0)=0.\end{eqnarray}
Let $\kappa=|\{\pi_i-\frac{1}{N}\neq 0:~i=1,2,...,N \}|$, i.e., the number of  elements in $\{\pi_i-\frac{1}{N}\neq 0:~i=1,2,...,N \}$. Clearly, $\kappa= 0,2,3,...,N$, by $\sum_{i=1}^N\pi_i=1$. $\kappa=0$ corresponds to Corollary \ref{corollary-1}. For $\kappa\geq 2$, without loss of generality, we assume that $\pi_i-\frac{1}{N}\neq 0, i=1,2,...,\kappa$. Thus, by (\ref{pi}), we have
\begin{eqnarray}
(\pi_1-\frac{1}{N})x_1(0)+\cdots+(\pi_{\kappa}-\frac{1}{N})x_{\kappa}(0)=0.\end{eqnarray}
Its solvable subspace is denoted by $V_1$, whose dimension is clearly
$\kappa-1$.  Therefore, for any $x(0)\in V_1 \oplus R^{N-\kappa}$,
the closed loop system (\ref{closed-loop-1}) will achieve the mean
square and strong average consensus.  It is worth pointing out
that the dimension of  $V_1\oplus R^{N-\kappa}$ is $N-1$.

\section{Networks with Switching Topology}\label{section-switching}

In this section, we extend results of last section to the case
that the network topology $G(t)$ of multi-agent systems is time-varying. The dependence of $G(t)$ on $t$ may be characterized by the switching signal $\sigma(t)$ in the meaning that if $\sigma(t)=k$, $G(t)=G_{\sigma(t)}=G^{(k)}$, $k=1,2,...T^*$. Here, the set of all possible digraphs is $\{G^{(k)}, k =
1,2,... ,T^*\}$. Therefore, the neighborhood of each node may vary with time. At any time $t\geq 0$, we can divide the nodes into two classes. On class is the isolated nodes denoted by $\mathcal{N}^s(t)$. The other class is denoted by $\mathcal{N}^a(t)$. For any node in $\mathcal{N}^a(t)$, either it has at least one neighbor, or it is a neighbor of node in $\mathcal{N}^a(t)$. In the following, we call nodes in $\mathcal{N}^a(t)$ to be active at time $t$. The active nodes with corresponding edges of $G(t)$ constructs a subgraph of $G(t)$, which is denoted by $G^a(t)$. While $G^s(t)$ denotes the graph composed by the isolated nodes of $G(t)$.  Therefore, digraph $G(t)$ may be viewed as the non-intersecting union of $G^a(t)$ and $G^s(t)$. Here, by the union of a collection
of graphs, we mean the graph whose nodes and edges set are the
unions of nodes and edge sets of the graphs in the collection. Similarly, we can define $G^{(k)a}$, $\mathcal{N}^{(k)a}$, and $G^{(k)s}$, $\mathcal{N}^{(k)s}$, $k=1,...,T^*$.
By (\ref{protoco0}) and (\ref{agent-dynamics}), through rearranging $x(t)$, we get the following closed loop system in
vector form
\begin{eqnarray}\label{closed-loop-switching-1}
\left\{\begin{array}{l}dx^a(t)=-aL^a(t)x^a(t)dt+a\Sigma^a(t)
y^a(t)dW^a(t),\\
\dot{x}^s(t)=0,\end{array}\right.\end{eqnarray} where, $x^a(t), L^a(t),
\alpha^a(t), \Sigma^a, y^a(t), W^a(t)$ correspond to subgraph
$G^a(t)$, and $x^s(t)$ corresponds to $G^s(t)$. Clearly, when $\sigma(t)=k$, the dimension of $x^a(t)$ is $|G^{(k)a}|$.

%If necessary, we reorder $x_i(t), i=1,2,\cdots, N$,
%such that the first $|\mathcal{N}^a(t)|$ elements in $x(t)$ being the
%states of the active agents, i.e., $x(t)=[x(t)^{aT}, x(t)^{sT}]^T$.
%Define $L(t)=\mbox{diag}(L^a(t), 0)$ with $0$ being
%$|\mathcal{N}(t)^s|\times|\mathcal{N}(t)^s|$ matrix. Therefore,
%\begin{eqnarray}
%dx(t)=-a(t)L(t)x(t)dt+\left(\begin{array}{c}a(t)\alpha^a(t)\Sigma^a
%y^a(t)\\0\end{array}\right)dW(t).
%\end{eqnarray}

To see the consensus of the agents, consider an infinite sequence of nonempty, bounded and contiguous
interval $[t_\tau,t_{\tau+1}), \tau=0,1,...$, starting at $t_0=0$ with
$t_{\tau+1}-t_\tau\leq T$, $T>0$. Suppose that during each interval
$[t_\tau,t_{\tau+1})$, $t_\tau^0(=t_\tau), t_\tau^1,..., t_\tau^{m_\tau}$ are the
switching instances of $\sigma(t)$, i.e., the points of
discontinuity of $\sigma(t)$, satisfying $\max\{t_\tau^{l+1}-t_\tau^l,
t_{\tau+1}-t_\tau^{m_\tau}\}\geq \bar{\tau}, \bar{\tau}>0, l=0,1,...,m_\tau-1$. To facilitate
the the following analysis, we relabel $t_{\tau+1}$ as $t_\tau^{m_\tau+1}$.
To this end, we need the following assumption.

(A2). For any $k=1,2,...$, the union of digraph $\{G(t), t_k\leq t<t_{k+1}\}$ contains a spanning tree.

Define the class of symmetric matrices:
$\mathcal{D}_N=\{D|D\in \mathbb{R}^{N\times N},D\geq 0, \mbox{Null}(D)=\mbox{span}\{1_N\}\}$.
Clearly, for any $D_1, D_2\in \mathcal{D}_N, \beta_1, \beta_2>0$,
$\beta_1D_1+\beta_2D_2\in\mathcal{D}_N$.
Notice the following properties.

(P1). For any $D\geq 0$, $x\in \mbox{Null}(D)$ if and only if $x^TDx=0$.

(P2). Let $D_1, D_2\in \mathcal{D}_N$, then there exists $c^*,c^{**}>0$ such that
$c^{*}D_2\leq D_1\leq c^{**}D_2$.

Clearly, (P1) can be derived by simple linear algebra knowledge. In fact, (P2) has been already presented implicitly in \cite{Huang-09-SIAM}\cite{Huang-10-Automatica}, and used in \cite{Huang-10-Automatica}. For the sake of completeness, a brief discussion about (P2) is given here. For any nonzero $x\in \mathbb{R}^N$, define $y=\frac{x}{|x|},\mathcal{B}_1=\{x:~|x|=1\}$, $\theta=\frac{\sqrt{N}}{N}1_N$, $\mathcal{B}_1^*=\mathcal{B}_1\backslash \{\theta\}$,  $\lambda_1=\max_{x\in \mathcal{B}_1^*}x^TD_1x$, $\lambda_2=\min_{x\in \mathcal{B}_1^*}x^TD_2x$. Then $x^TD_1x=|x|^2y^TD_1y\leq \lambda_1 |x|^2\leq \frac{\lambda_1}{\lambda_2}\lambda_2|x|^2\leq \frac{\lambda_1}{\lambda_2}x^TD_2x$. Therefore, we may select $c^{**}$ as $\frac{\lambda_1}{\lambda_2}$. $c^{*}$ can be similarly constructed.

Clearly, the same edge $(j,i)$ may be present in some of $G^{(k)},
k=1,2,...,T^*$. Denote the number of  $G^{(k)}$ that the edge $(j,i)$ is present in $G^{(k)},
k=1,2,...,T^*$, by $e_{ji}$, and let $e=\max_{i,j}e_{ji}$. Clearly, $e\geq 1$.  Then we have the following result.

\begin{lemma}\label{lemma-huang-2}
Let $U=\frac{1}{2}\sum_{i=1}^N\sum_{j=1}^N|x_j-x_i|^2$,
$P^{(k)}=\frac{1}{2}\sum_{i\in \mathcal{N}^{(k)a}}\sum_{j\in
\mathcal{N}^{(k)a}_{i}}|x_j-x_i|^2, k=1,2,...,T^*$. Then the union of $\{G^{(k)},
k=1,2,...,T^*\}$ contains a spanning tree if and only if the following
is satisfied
\begin{eqnarray}\label{lemma-huang-2-0}
c_{*}U\leq\sum_{k=1}^{T^*}P^{(k)},
\end{eqnarray}
where $c_{*}>0$. In addition, the ``only if " part can be strengthened to
\begin{eqnarray}\label{lemma-huang-2-0-0}
c^{*}U\leq\sum_{k=1}^{T^*}P^{(k)}\leq ec^{**}U,
\end{eqnarray}
where $c^*, c^{**}$ are given in (P2).
\end{lemma}

\emph{{Proof}}. \emph{Necessity.} Denote the union of $\{G^{(k)},
k=1,2,...,T^*\}$  by $G^u$. A simple calculation shows that
\begin{eqnarray}\label{P}
\frac{1}{e}\sum_{k=1}^{T^*}P^{(k)} \leq P^u \equiv \frac{1}{2} \sum_{i\in\mathcal{N}(G^u)}\sum_{j\in\mathcal{N}(G^u)_i}|x_j-x_i|^2\leq  \sum_{k=1}^{T^*}P^{(k)},
\end{eqnarray}
where $\mathcal{N}(G^u)$ denotes the nodes set of $G^u$, $\mathcal{N}(G^u)_i$ is the neighborhood of node $i$ in $G^u$. Because $G^u$ contains a spanning tree, every one in $\{x_i, i=1,2,...,N\}$ will appear in $P^u$. By the special structure of $P^u$, we have that $P^u\equiv 0$ if and only if $[x_1,...,x_N]^T=x\in \mbox{span}\{1_N\}$. On the other hand, as $P^u$ is a quadratic form, there exists $H\geq 0$ such that  $P^u=x^THx$. By (P1), $\mbox{Null}(H)=\mbox{span}\{1_N\}$, and thus $H\in \mathcal{D}_N$. As $U=x^TSx$ with $S=NI_N-1_N1_N^T\in \mathcal{D}_N$, by property (P2) and (\ref{P}), (\ref{lemma-huang-2-0}) and (\ref{lemma-huang-2-0-0}) are followed with $c_{*}=c^{*}$.

\emph{Sufficiency.} Define $H$ such that $P^u=x^THx$. Clearly, $H\geq 0$. By (\ref{lemma-huang-2-0}), $\sum_{k=1}^{T^*}P^{(k)}=0$ implies that $U= 0$, and thus $x\in\mbox{span}\{1_N\}$. By property (P1), $H$ must be in $\mathcal{D}_N$. Now, we show that $G^u$ contains a spanning tree. This is proved by contradiction. Assume that $G^u$  does not contain a spanning tree.
Without loss of generality, suppose that $G^u$ contains two non-intersecting subgraphs  $G^{u1}$
and $G^{u2}$, both of which contain a spanning tree.
Therefore,
\begin{eqnarray}\label{lemma-huang-0-inequ-2}
P^{u}=x^THx=x^{(1)T}H_1x^{(1)}+x^{(2)T}H_2x^{(2)},
\end{eqnarray}
where $\{x^{(1)},H_1\}$, and $\{x^{(2)}, H_2\}$ correspond to $G^{u1}$ and $G^{u2}$, respectively; in addition, $x^{(1)T}=(x^{(1)})^T$, $x^{(2)T}=(x^{(2)})^T$. Selecting $x^{(1)}=c_{*1}1_{|\mathcal{N}(G')|}$, $x^{(1)}=c_{*2}1_{|\mathcal{N}(G'')|}$ with $c_{*1}\neq c_{*2}$, we can assert that the right side of (\ref{lemma-huang-0-inequ-2}) is equal to zero. This contradicts that $\mbox{Null}(H)=\mbox{span}(1_N)$ by (P1). Therefore, $G^u$ must contain a spanning tree. This completes the proof. \hfill$\square$

\begin{theorem}\label{theorem-switching-mean-square}
Under assumption (A2), if for every $k=1,2,...,T^*$, $G^{(k)}$ is balanced, and $0<a<\bar{\bar{a}}$ with $\bar{\bar{a}}=\frac{N}{(N-1)\max_{i,j}\sigma_{ij}^2}$, then the distributed protocol (\ref{protoco0}) is a mean square and strong consensus protocol.
\end{theorem}

\emph{{Proof}}. Let $U(t)=\frac{1}{2}\sum_{i=1}^N\sum_{j=1}^N|x_j(t)-x_i(t)|^2 $. Clearly, $U(t)=x^T(t) S x(t),$ with $S=NI_N-1_N1_N^T$.
%\begin{eqnarray}\label{P(t)}
%U(t)&=&\frac{1}{2}\sum_{i=1}^N\sum_{j=1}^N\left[(x_j(t))^2+(x_i(t))^2-2x_j(t)x_i(t)\right]
%=x^T(t) S x(t),
%\end{eqnarray}
%where $S=NI_N-1_N1_N^T$.
As $G^{(k)}, k=1,2,...,T^*$, are balanced, we have
\begin{eqnarray*}
P(t) \equiv \frac{1}{2}\sum_{i\in \mathcal{N}^a(t)}\sum_{j\in \mathcal{N}_{ti}^a}|x_j(t)-x_i(t)|^2
=\frac{1}{2} x^{aT}(t)(D^a(t)+\bar{D}^a(t)-2A^a(t))x^a(t)=x^{aT}(t)L^a(t)x^a(t),
\end{eqnarray*}
where $x(t)^{aT}$ denotes $(x^a(t))^T$, $D^a(t)$ and $\bar{D}^a(t)$ are the in-degree and out-degree matrices of $G^a(t)$, and thus of $G(t)$.
As $U(t)=[(x^{a}(t))^T~~(x^s(t))^T]S[(x^{a}(t))^T~~(x^s(t))^T]^T$, by It\^o's formula, we have
\begin{eqnarray}\label{switching-P}
dU(t)&=&\left[-2Nax(t)^{aT}L^a(t)x^a(t)+a^2Tr\left(S^a(t)(\Sigma^a(t)
y^a(t))(\Sigma^a(t)
y^a(t))^T\right)\right]dt\nonumber \\
&&+2ax(t)^{aT}S^a(t)\Sigma^a(t) y^a(t)dW^a(t),
\end{eqnarray}
where
$S^a(t)=NI_{|\mathcal{N}^a(t)|}-1_{|\mathcal{N}^a(t)|}1^T_{|\mathcal{N}^a(t)|}\geq 0$,
and the properties that $1_{|\mathcal{N}^a(t)|}^TL^a(t)=0,
L^a(t)1_{|\mathcal{N}^a(t)|}=0$ are used for several times. Clearly, $y^a(t)$ is composed of elements with form $|x_j(t)-x_i(t)|$, where agent $j$ is a neighbor of agent $i$ at time $t$. A simple
calculation shows that
\begin{eqnarray}\label{switching-inequality-1}
&&Tr\left(S^a(t)(\Sigma^a(t)
y^a(t))(\Sigma^a(t) y^a(t))^T\right)=\sum_{i\in\mathcal{N}^a(t)}[S^a(t)]_{ii}\sum_{j\in\mathcal{N}^a_{ti}} \left(\sigma_{ji}|x_j(t)-x_i(t)|\right)^2\nonumber\\
&&\leq 2(N-1)\max_{i,j}\sigma_{ij}^2\cdot \frac{1}{2}\sum_{i\in\mathcal{N}^a(t)}\sum_{j\in\mathcal{N}^a_{ti}} |x_j(t)-x_i(t)|^2\nonumber\\
&&\equiv c_4 P(t), \end{eqnarray}
where $[S^a(t)]_{ii}=N-1$. Therefore,
\begin{eqnarray}\label{switching-U}
U(t)\leq U(0)-\int_0^t(2Na-c_4a^2)P(s)ds+\int_0^t2a(x^a(s))^TS^a(s)\Sigma^a(s)
y^a(s)dW^a(s).
\end{eqnarray}
Select $a$ such that
\begin{eqnarray*}
\gamma_2\equiv 2Na-c_4a^2>0,
\end{eqnarray*}
i.e.,
\begin{eqnarray}\label{a-switching-1}
0<a <\frac{2N}{c_4}=\bar{\bar{a}}.
\end{eqnarray}
By (\ref{martingle-1}), $\int_0^t2a(x^a(s))^TS^a(s)\Sigma^a(s)y(s)^adW^a(s)$ is a martingale. {By similar analysis to that of Theorem \ref{theorem-1}, we have that $EU(t)$ is {differentiable a.e. with respect to $t$}}, and
\begin{eqnarray}\label{U-2}
\frac{d EU(t)}{dt}\leq -\gamma_2 EP(t)\leq 0.
\end{eqnarray}
Therefore, $EU(t)$ converges, and denote the limitation of $EU(t)$ by $\bar{U}$.
%\begin{eqnarray}\label{u-3}
%\lim_{t\rightarrow \infty}EU(t)=\bar{U}\geq 0.
%\end{eqnarray}
We will show that $\bar{U}=0$. And this is proved by
contradiction. Assume that $\bar{U}>0$.  Define
\begin{eqnarray*}
P^{(k)}(t)=\frac{1}{2}\sum_{i\in \mathcal{N}^{(k)a}}\sum_{j\in
\mathcal{N}^{(k)a}_{i}}|x_j(t)-x_i(t)|^2,~ k=1,2,...,T^*.\end{eqnarray*}
By assumption (A2) and Lemma \ref{lemma-huang-2}, we have
\begin{eqnarray*}
\sum_{k=1}^{T^*}EP^{(k)}(t)\geq c_* EU(t),~~c_*>0.\end{eqnarray*}
Therefore,
\begin{eqnarray*}
\liminf_{t\rightarrow \infty}\sum_{k=1}^{T^*}EP^{(k)}(t)\geq
c_*\bar{U}>0,\end{eqnarray*} which implies that for some $k^*\in
\{1,2,...,T^*\}$,
\begin{eqnarray*}
\liminf_{t\rightarrow \infty}EP^{(k^*)}(t)\geq
\frac{c_*\bar{U}}{T^*}.\end{eqnarray*} Thus, there exists
$T_1>0$ such that
\begin{eqnarray}\label{U-3}
EP^{(k^*)}(t)\geq \frac{c_*\bar{U}}{2T^*},~~ t\geq
T_1.\end{eqnarray}
By assumption (A2), we know that there exists a subinterval $[t^l_\tau,t^{l+1}_{\tau+1})$ of $[t_\tau,t_{\tau+1}), \tau=1,2,...$, during which $G(t)\equiv G^{k^*}$. On this subinterval $[t^l_\tau,t^{l+1}_{\tau+1})$, by (\ref{U-2}), it follows that
\begin{eqnarray}\label{U-4}
\frac{d EU(t)}{dt}\leq -\gamma_2 EP(t)^{k^*}\leq -\gamma_2\frac{c_*\bar{U}}{2T^*}.
\end{eqnarray}
Clearly, there exists $\tau_{T_1}$, such that $t_{\tau_{T_1}}>T_1$. Combining (\ref{U-3})(\ref{U-4}), we have that for $t>t_{\tau_{T_1}}$,
\begin{eqnarray}\label{U-10}
EU(t)\leq EU_{t_\tau}\leq EU_{t_{\tau-1}}-\bar{\tau}\cdot \gamma_2
\frac{c_*\bar{U}}{2T^*}\leq\cdots\leq EU_{t_{\tau_{T_1}}}-(\tau-\tau_{T_1})\bar{\tau}\cdot
\gamma_2 \frac{c_*\bar{U}}{2T^*}, \end{eqnarray} where $\tau$ is the
largest integer such that $t\geq t_\tau$. In (\ref{U-10}), the
inequality
\begin{eqnarray*}
EU_{t_\tau}\leq EU_{t_{\tau-1}}-\int_{t_{\tau-1}}^{t_\tau}\gamma_2 EP^*(s)ds\leq
EU_{t_{\tau-1}}-\gamma_2 \bar{\tau} \frac{c_*\bar{U}}{2T^*}\end{eqnarray*} is
used iteratively. Letting $t\rightarrow \infty$ in (\ref{U-10}), we
have that $\lim_{t\rightarrow \infty}EU(t)=-\infty$. This contradicts
that $EU(t)\geq 0$. Therefore, we have that
\begin{eqnarray}
\lim_{t\rightarrow \infty}EU(t)=0,
\end{eqnarray}
i.e.,
\begin{eqnarray*}
\lim_{t\rightarrow\infty}E|x_j(t)-x_i(t)|^2=0, ~~i, j=1,2,...,N.
\end{eqnarray*}
Notice that
\begin{eqnarray*}
d(1_N^Tx(t))=a1_{|\mathcal{N}^a(t)|}^T\Sigma^a(t)y^a(t)dW^a(t),
\end{eqnarray*}
i.e.,
\begin{eqnarray}\label{1x}1_{N}^Tx(t)=1_N^Tx(0)+\int_0^ta1_{|\mathcal{N}(s)^a|}^T\Sigma^a(s)y^a(s)dW^a(s).\end{eqnarray}
Similar to (\ref{switching-inequality-1}), we have that
\begin{eqnarray}\label{switching-inequality-2}
(1_{|\mathcal{N}(s)^a|}^T\Sigma^a(s)y^a(s))^T1_{|\mathcal{N}(s)^a|}^T\Sigma^a(s)y^a(s)
\leq c_5 P(s),\end{eqnarray} for some $c_5>0$.  On the other hand, by
(\ref{U-2}), it follows that
\begin{eqnarray}\label{switching-inequality-3}
\gamma_2\int_0^\infty EP(s)ds\leq EU(0)-EU(\infty),
\end{eqnarray}
where $EU(\infty) \equiv \lim_{t\rightarrow \infty}EU(t)=0$. Therefore,
\begin{eqnarray*}
\sup_{t\geq 0}E|1^T_Nx(t)|^2<\infty.
\end{eqnarray*}
By similar analysis to that of Theorem \ref{theorem-1}, $1^T_Nx(t)$ converges to $1^Tx(\infty)$ almost surely and in the sense of mean square, where $1_N^Tx(\infty)$ is defined as
\begin{eqnarray*}
1_{N}^Tx(\infty)=1_N^Tx(0)+\int_0^\infty
a1_{|\mathcal{N}(s)^a|}^T\Sigma^a(s)y^a(s)dW^a(s).
\end{eqnarray*}
Therefore, by the fact
\begin{eqnarray*}
|x_i(t)-\frac{1}{N}1_N^Tx(\infty)|^2\leq 2|x_i(t)-\frac{1}{N}1_N^Tx(t)|^2+2| \frac{1}{N}1_N^Tx(t)-\frac{1}{N}1_N^Tx(\infty)|^2,
\end{eqnarray*}
we can conclude that $ x_i(t)$ converges to $\frac{1}{N}1^T_Nx(\infty)$ in the sense of mean square, for any $i=1,2,...,N$.

Similar to Theorem
\ref{theorem-2}, we can easily prove that (\ref{protoco0}) is a strong consensus protocol. And the proof is omitted here. \hfill$\square$

\begin{corollary}
Under the conditions of Theorem \ref{theorem-switching-mean-square},
the protocol (\ref{protoco0}) is an asymptotic unbiased mean
square (strong, respectively) average consensus protocol.
\end{corollary}

\emph{{Proof}}. Clearly,
\begin{eqnarray*}
E\frac{1}{N}1^T_Nx(\infty)=\frac{1}{N}1_N^Tx(0),
\end{eqnarray*}
and the variance is
\begin{eqnarray*}
Var(\frac{1}{N}1^T_Nx(\infty))=E\int_0^\infty(1_{|\mathcal{N}(s)^a|}^T\Sigma^a(s)y^a(s))^T
1_{|\mathcal{N}(s)^a|}^T\Sigma^a(s)y^a(s)ds<\infty.
\end{eqnarray*}
The conclusion follows easily.
\hfill $\square$

\begin{remark}If the noises disappear during communication processes, i.e., $\sigma_{ji}=0, i,j=1,...,N,$ in
(\ref{measurement-0}), then under protocol (\ref{protoco0}) the
closed-loop system is
\begin{eqnarray}\label{deterministic-switch}
\dot{x}(t)=-aL(t)x(t),
\end{eqnarray}
which is the objective of \cite{Jadbabaie}, \cite{Moreau}, \cite{Ren
Wei}, \cite{Wang-long}.
Comparing to above Theorem \ref{theorem-switching-mean-square}, results in \cite{Moreau}, \cite{Ren
Wei}, \cite{Wang-long} do not need the balancedness condition.  It should be
mentioned that existing results are all based on the fact: the explicit solution  of  (\ref{deterministic-switch}) can be easily expressed.  While for the noise-driven case, it is always impossible to derive the explicit solution of the closed-loop stochastic differential equation. Therefore, we adopt the Lyapunov-based approach to tackle this problem. As a cost, balanceness property of the digraphs is needed. {For more about Lyapunov-based approach for consensus of  agents, readers may refer to \cite{Yu} and references therein. In \cite{Yu}, the authors study the consensus problems in direct networks with nonlinear dynamics.}
\end{remark}

\section{Simulation}
\label{simulation}

%\subsection{Simulations for Fixed topology}

\textbf{Example 5.1} Consider a dynamic network of four agents with
fixed topology with  $\mathcal{N}=\{1,2,3,4\}$,
$\mathcal{E}=\{(1,3),(2,1),(3,2),(3,1),(4,3)\}$. The quotient
digraph is shown in Fig.1. $\sigma_{31}=\sigma_{12}=\sigma_{13}=\sigma_{23}=\sigma_{34}=1$. The consensus gain $a$ is selected as 0.05. For initial states $x^1(0)=1, x(0)^2=20,
x(0)^3=50,x(0)^4=-5$, under protocol (\ref{protoco0}), the states
of the closed loop system are shown in Fig.2.
From Fig.2, we can see that when $t$ sufficiently large, the
common value of the sample path of the agents' states is about
21.21; the corresponding mean value is about 17.98; while the average
of the initial states of all agents is 16.5. This means that the
multi-agent system does not achieve (global) average consensus.
However, local average consensus may be attained for some
the initial values of the states of the agents. By simple
computation, the $\pi$ defined in (\ref{pi}) is $(0.5,0.25,0.25,0)$.
Therefore, if $x_1(0)$ is equal to $x_4(0)$, the average
consensus can be achieved. This is validated by Fig.3 with initial
states $x_1(0)=-5, x_2(0)=20, x_3(0)=50,x_4(0)=-5$.

\textbf{Example 5.2} Consider a dynamic network of four agents with
undirected network topology. The network topology are changed as
follows: when $t=[2k,2k+1), k=0,1,2...$, it has the structure
shown in Fig.4.(a); while when $t=[2k+1,2k+2), k=0,1,2,...$, it is
described  by (b) of Fig.4. Assume that all the
$\sigma_{..}$ are equal to 1. $a$ defined in (\ref{a-switching-1}) is
selected to be 2. The initial state $x(0)$ is $[1~  2~  5~-10]^T$. When
(\ref{protoco0}) is applied, the closed loop states will reach a
consensus shown by Fig.5.

\setlength{\unitlength}{1cm}
\begin{center}\begin{picture}(3,3)
\put(0.45,0.45){\circle{0.45}\makebox(-0.9,0){4}}
\put(2.55,0.45){\circle{0.45}\makebox(-0.9,0){3}}
\put(0.45,2.55){\circle{0.45}\makebox(-0.9,0){1}}
\put(2.55,2.55){\circle{0.45}\makebox(-0.9,0){2}}
\put(2.325,0.45){\vector(-1,0){1.65}}
\put(2.55,2.325){\vector(0,-1){1.65}}
\put(0.675,2.55){\vector(1,0){1.65}}
\put(2.3910,0.6090){\vector(-1,1){1.77}}
\end{picture}\\
Fig.1. The fixed topology of Example 5.1
\end{center}

\begin{figure}[h]
 \centering{ \includegraphics[height=0.35\textwidth]{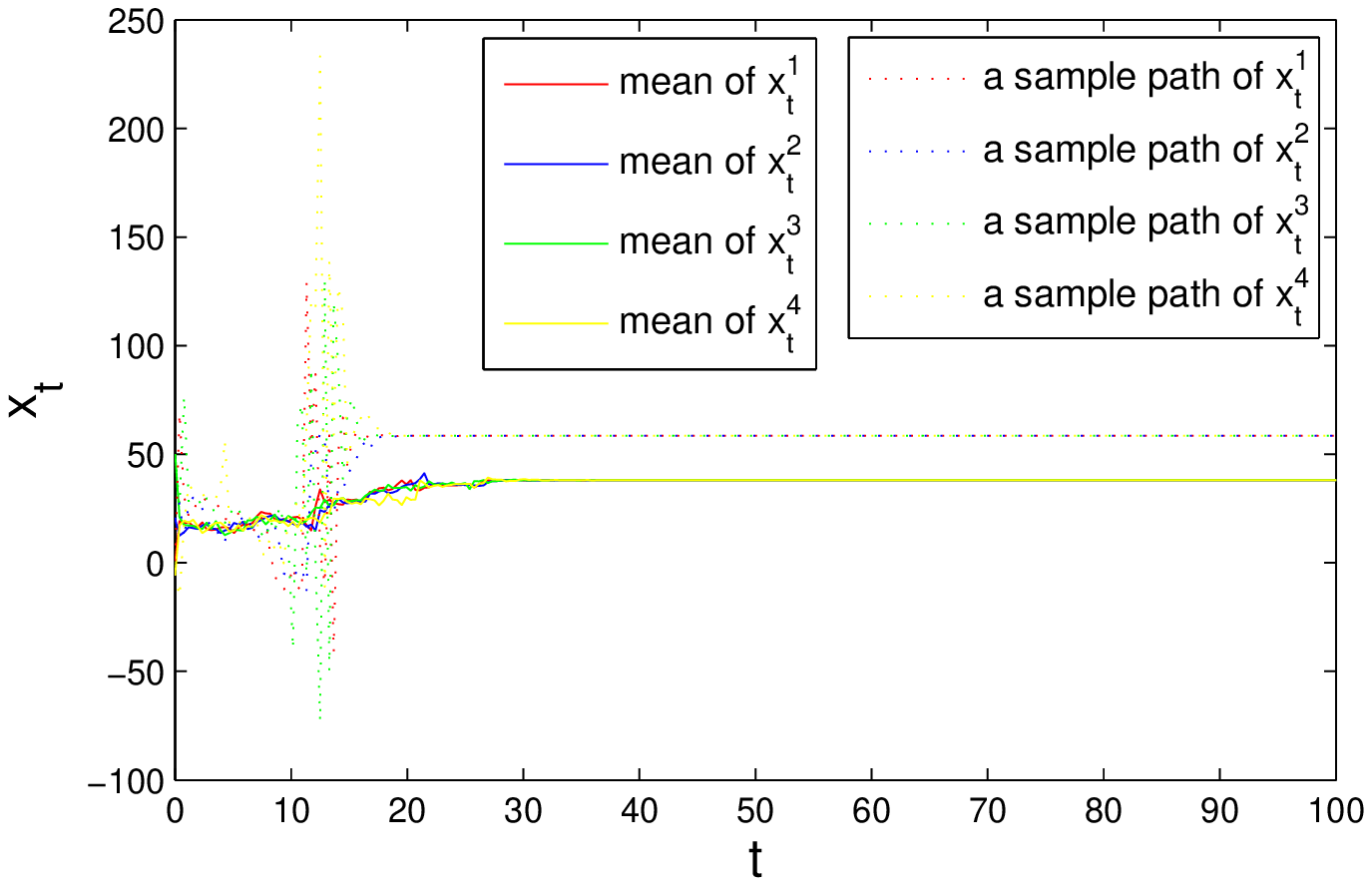}}\\
\begin{center}Fig.2. Curves of the states of Example 5.1 \\
with initial state [1, 20, 50, -5]
\end{center}
\end{figure}
\begin{figure}[h]
\centering{ \includegraphics[height=0.35\textwidth]{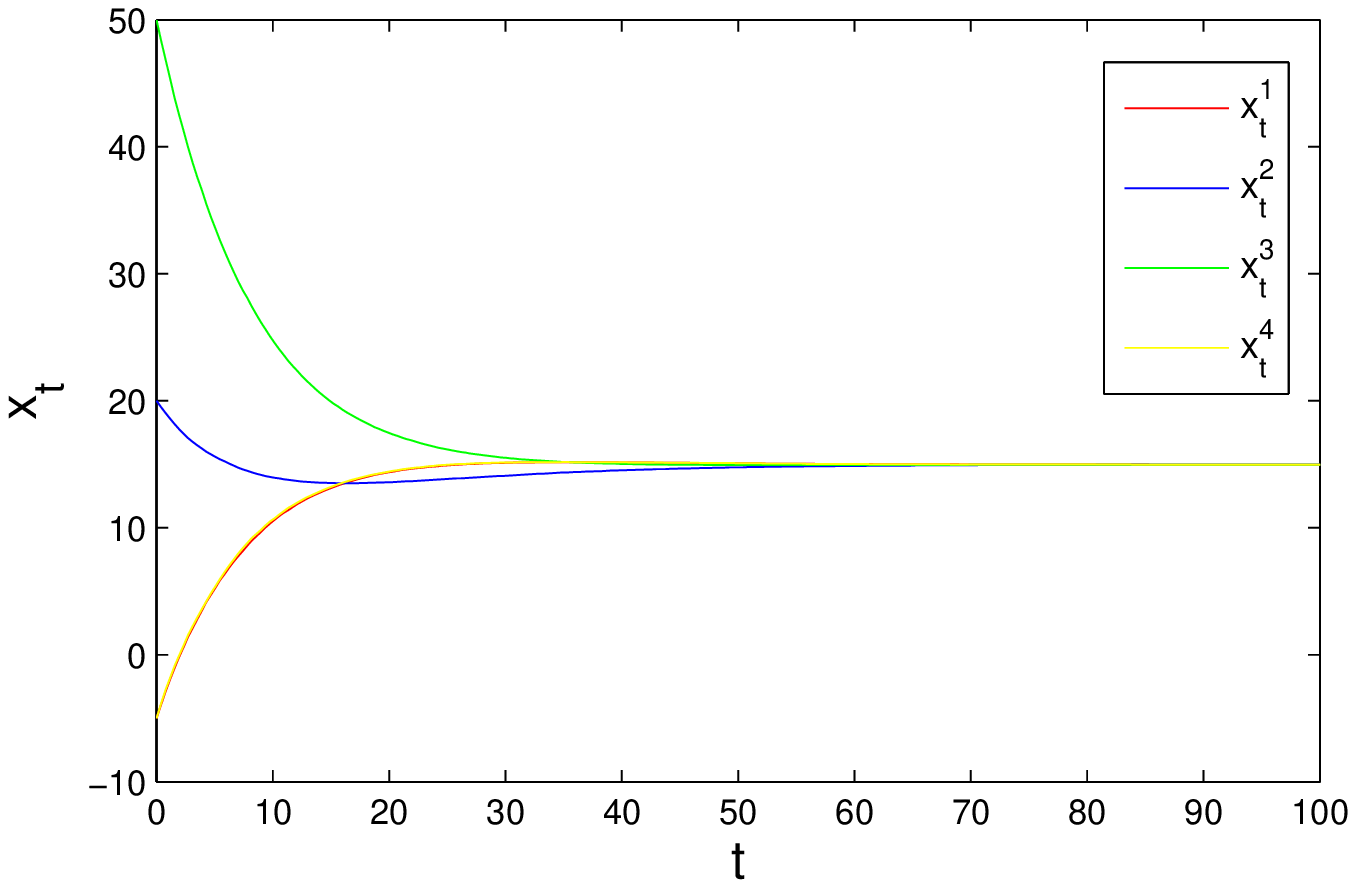}}\\
\begin{center}Fig.3. Curves of the states' means of Example 5.1\\
with initial state [-5, 20, 50, -5]
\end{center}
\end{figure}

\setlength{\unitlength}{1cm}
\begin{center}\begin{picture}(7,3)
\put(0.45,0.45){\circle{0.45}\makebox(-0.9,0){4}}
\put(2.55,0.45){\circle{0.45}\makebox(-0.9,0){3}}
\put(0.45,2.55){\circle{0.45}\makebox(-0.9,0){1}}
\put(2.55,2.55){\circle{0.45}\makebox(-0.9,0){2}}
%\put(2.6,0.6){\vector(-1,0){1.7}}
\put(2.55,2.325){\line(0,-1){1.65}}
\put(0.675,2.55){\line(1,0){1.65}}
%\put(2.6879,0.8121){\vector(-1,1){1.885}}

\put(4.45,0.45){\circle{0.45}\makebox(-0.9,0){4}}
\put(6.55,0.45){\circle{0.45}\makebox(-0.9,0){3}}
\put(4.45,2.55){\circle{0.45}\makebox(-0.9,0){1}}
\put(6.55,2.55){\circle{0.45}\makebox(-0.9,0){2}}
\put(6.325,0.45){\line(-1,0){1.65}} %\put(6.9,2.6){\vector(0,-1){1.7}}
%\put(4.9,2.9){\vector(1,0){1.7}}
\put(6.3910,0.6090){\line(-1,1){1.77}}
\end{picture}\\
(a)~~~~~~~~~~~~~~~~~~~~~~~~(b)\\\vspace{0.5em} Fig.4.
The switching topology of Example 5.2
\end{center}

\begin{figure}[h]
 \centering{ \includegraphics[height=0.35\textwidth]{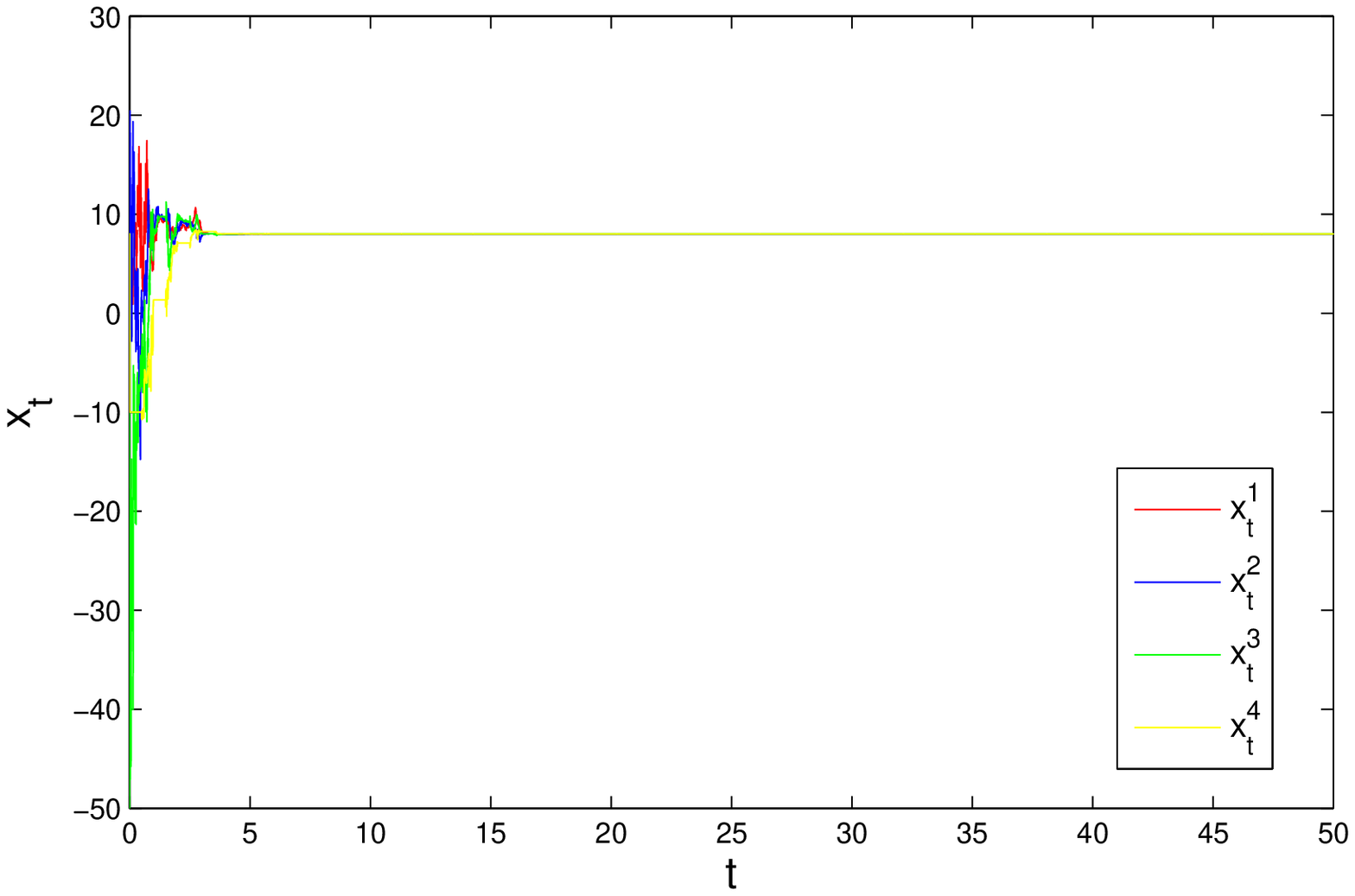}}\\
\centerline{Fig.5. Curves of the states of Example 5.2}
\end{figure}

\section{conclusion}
\label{conclusion}

This paper considers the consensus problem of first order integrator systems under uncertainty environment. The measurement noises are modeled to be multiplicative. For fixed and switching topologies cases, mean square and strong consensus are achieved. For further research, consensus problems under measurement noises with leaders are valuable for some applied scenarios. In addition, it is an issue to consider the stochastic varying network topology cases.

%\bibliographystyle{plain}        % Include this if you use bibtex
%\bibliography{autosam}           % and a bib file to produce the
                                 % bibliography (preferred). The
                                 % correct style is generated by
                                 % Elsevier at the time of printing.

%% main text
%\section{}
%\label{}

%% The Appendices part is started with the command \appendix;
%% appendix sections are then done as normal sections
%% \appendix

%% \section{}
%% \label{}

%% References
%%
%% Following citation commands can be used in the body text:
%% Usage of \cite is as follows:
%%   \cite{key}          ==>>  [#]
%%   \cite[chap. 2]{key} ==>>  [#, chap. 2]
%%   \citet{key}         ==>>  Author [#]

%% References with bibTeX database:

%\bibliographystyle{model1-num-names}
%\bibliography{<your-bib-database>}

%% Authors are advised to submit their bibtex database files. They are
%% requested to list a bibtex style file in the manuscript if they do
%% not want to use model1-num-names.bst.

%% References without bibTeX database:

% \begin{thebibliography}{00}

%% \bibitem must have the following form:
%%   \bibitem{key}...
%%

% \bibitem{}

% \end{thebibliography}

\end{document}